\documentclass[12pt]{amsart}
\usepackage{pifont}
\usepackage{mathrsfs}
\usepackage{geometry}
\usepackage{titletoc}
\usepackage{stix2}

\usepackage{amsmath}
\usepackage{amssymb} 
\usepackage{enumitem} 
\usepackage{mathtools} 
\usepackage[table]{xcolor} 
\usepackage[all]{xy} 
\usepackage{tikz} 
\usepackage{tikz-cd}
\usepackage{indentfirst} 
\usepackage{babel} 
\usepackage{setspace} 

\usepackage[colorlinks,linkcolor=red,anchorcolor=green,citecolor=blue]{hyperref} 
\hypersetup{linktocpage = true} 

\usepackage{rotating} 

\usepackage{ytableau} 
\usepackage{longtable} 
\newcolumntype{M}[1]{>{\centering\arraybackslash}m{#1}} 

\usepackage{fancyhdr}

\newcommand\bp{{\bar\partial}}

\theoremstyle{plain}
\newtheorem{thm}{Theorem}[section]
\newtheorem{lemma}[thm]{Lemma}
\newtheorem{prop}[thm]{Proposition}
\newtheorem{cor}[thm]{Corollary}
\newtheorem{defn}[thm]{Definition}

\theoremstyle{definition}
\newtheorem{example}[thm]{Example}
\newtheorem{remark}[thm]{Remark}

\newcommand{\btheorem}{\begin{thm}}
	\newcommand{\etheorem}{\end{thm}}
\newcommand{\bproposition}{\begin{prop}}
	\newcommand{\eproposition}{\end{prop}}
\newcommand{\bdefinition}{\begin{defn}}
	\newcommand{\edefinition}{\end{defn}}
\newcommand{\bcorollary}{\begin{cor}}
	\newcommand{\ecorollary}{\end{cor}}
\newcommand{\bproof}{\begin{proof}}
	\newcommand{\eproof}{\end{proof}}
\newcommand{\bremark}{\begin{remark}}
	\newcommand{\eremark}{\end{remark}}
\newcommand{\eexample}{\end{example}}
\newcommand{\bexample}{\begin{example}}

\newcommand{\elemma}{\end{lemma}}
\newcommand{\blemma}{\begin{lemma}}

\newcommand{\sq}{\sqrt{-1}}
\newcommand{\suml}{\sum\limits}
\newcommand{\p}{\partial}

\renewcommand{\bar}{\overline}

\renewcommand{\phi}{\varphi}

\newcommand{\beq}{\begin{equation}}
\newcommand{\eeq}{\end{equation}}
\newcommand{\ee}{\end{eqnarray*}}
\newcommand{\be}{\begin{eqnarray*}}

\newcommand{\bd}{\begin{enumerate}}
	\newcommand{\ed}{\end{enumerate}}

\renewcommand{\hat}{\widehat}

\newcommand{\qtq}[1]{\quad\mbox{#1}\quad}
\renewcommand{\bp}{\bar{\partial}}

\newcommand{\ts}{\otimes}

\renewcommand{\S}{{\mathbb S}}

\newcommand{\Z}{{\mathbb Z}}
\renewcommand{\>}{\rightarrow}

\newcommand{\ssR}{{\mathfrak R}}

\newcommand{\C}{{\mathbb C}}

\newcommand{\R}{{\mathbb R}}

\setlist[itemize]{leftmargin=*}
\setlist[enumerate]{leftmargin=*}

\numberwithin{equation}{section} 

\makeatletter

\usepackage{fancyhdr}
\title{Manifolds with non-positive second Chern-Ricci curvature}

\author{Xiaokui Yang}

\address{Xiaokui Yang, Department of Mathematics and Yau Mathematical Sciences Center, Tsinghua University, Beijing, 100084, China}
\email{xkyang@mail.tsinghua.edu.cn}

\begin{document}
	
	\begin{abstract} In this paper, we establish a K\"ahlerian or projectivity criterion for a class of compact Hermitian surfaces with non-positive second Chern-Ricci curvature.
	\end{abstract}

	\maketitle

		\section{Introduction} 
	
Let $(M,\omega)$ be a compact Hermitian manifold. When $\omega$ is not K\"ahler, there are many different connections on $T^{1,0}M$ and their curvature tensors are different. Recently, there are many works on curvature relations for various connections, for instances, \cite{Broder2023}, \cite{He2020}, \cite{LiuYang2012}, \cite{LiuYang2017}, \cite{Ni2023},  \cite{WangYang2019} and \cite{YangZheng2018}.
 It is still a challenge to investigate the geometry and topology determined by those connections and curvature tensors.  The key difficulty is that the metric differential has no K\"ahler symmetry. Hermitian metrics with special symmetric are investigated extensively in recent years.  For instances, certain balanced metrics are constructed  in  \cite{Fu2012a},  and the K\"ahlerian of various K\"ahler-like metrics are obtained in \cite{Angella2022}, \cite{Podesta2018}, \cite{Podesta2023}, \cite{Yau2023}, \cite{Zhao2023}. For more discussions  along this line,  one can see  \cite{Correa2023}, \cite{Fino2004}, \cite{Fu2010},  \cite{Fu2012}, \cite{Ivanov2013} and  \cite{Zheng2019}.\\ 
 
 On the other hand, the geometry and topology of compact Hermitian surfaces are well studied  by using the  Weyl structure and conformal transformations. In particular, when $(M,\omega)$ is Hermitian-Einstein and $\omega$ is a Gauduchon metric, Gauduchon and Ivanov observed in \cite{Gauduchon1997} (see also \cite{Tod1992}) that the Lee form is parallel and such manifolds are classified  in \cite{Gauduchon1995} and \cite{Gauduchon1997}. For more results along this line, we refer to
 \cite{Apostolov1996}, \cite{Apostolov}, \cite{Calderbank1999},  \cite{LeBrun1999}, \cite{LeBrun2012}, \cite{Pedersen1993} and the references therein.\\


  It is well-known that  the first Chern-Ricci curvature $\Theta^{(1)}$ of the  Chern connection represents the first Chern class of  $M$.  In particular, if $\omega$ has positive first Chern-Ricci curvature  $\Theta^{(1)}$, then $M$ must be a K\"ahler manifold with  $c_1(M)>0$. Moreover, it has a K\"ahler metric with positive Ricci curvature (\cite{Yau1978}) and so $M$ is simply connected. However, the geometry and topology of the second Chern-Ricci curvature $\Theta^{(2)}$ is more mysterious.  For instance,
 it is well-known that (e.g. \cite{LiuYang2012, LiuYang2017}) the non-K\"ahler manifold $\S^3\times \S^1$ has positive  $\Theta^{(2)}$ and $\pi_1(M)\cong \Z$. On the other hand, we proved in \cite{Yang2018} that if a compact K\"ahler manifold admits a Hermitian metric with positive second Chern-Ricci curvature $\Theta^{(2)}$, then $M$ is projective and simply connected.\\

   In this paper, we establish a K\"ahlerian or projectivity criterion for a class of compact Hermitian manifolds with non-positive second Chern-Ricci curvature. We begin with the classical result of Gauduchon-Ivanov \cite{Gauduchon1997}.

	\btheorem\label{main0} Let $(M,\omega)$ be a compact Hermitian surface. Suppose $\omega$ is  Gauduchon and weakly Hermitian-Einstein, i.e. $\p\bp\omega=0$ and
	\beq \Theta^{(2)}(\omega)=u \omega \eeq 
	for some $u\in C^\infty(M,\R)$. If $u\leq 0$, then $\omega$ is a K\"ahler metric.
	\etheorem

	\noindent The original proof of this theorem uses the geometry of Weyl-Einstein structure. 
	We shall give a simplified proof of this result and derive some applications.  For higher dimensional Hermitian manifolds, we obtain

	\btheorem\label{secondChernRicci} Let $M$ be a compact complex manifold with $\dim M=n$. If it admits a Hermtian metric $\omega$ with $\Lambda_\omega\p\bp\omega=0$ and 
	$\mathrm{\Theta^{(2)}}(\omega)$ is quasi-negative, then
	the signed top intersection number $$(-1)^n c_1^n(M)>0.$$  Moreover,  the cohomology groups $H^2_{\mathrm{dR}}(M,\R)$, $H^{1,1}_{\bp}(M,\C)$ and $H^{1,1}_{\mathrm{BC}}(M,\C)$  are all non-zero. 
	\etheorem

\noindent The key ingredient in the proof of Theorem \ref{secondChernRicci} is a combinatoric version of the Chern-Weil theory (Theorem \ref{thm:ChernWeil}).  As an application of Theorem \ref{secondChernRicci}, we obtain an extension of Theorem \ref{main0}:

	\bcorollary\label{surface}  Let $(M,\omega)$ be a compact Hermitian surface. If  $\omega$ is Gauduchon and  $$\mathrm{\Theta^{(2)}}(\omega)\leq 0,$$ then either  $(M,\omega)$ is a K\"ahler manifold with $c_1^2(M)=0$,  or $M$ is a projective manifold with $c_1^2(M)>0$.
	\ecorollary

	\noindent It worths to point out that similar results as in Theorem \ref{secondChernRicci} and Corollary \ref{surface} can not hold for compact Hermitian manifolds with $\mathrm{\Theta^{(2)}}>0$.  Indeed,  it is shown in \cite[Proposition~6.1]{LiuYang2012} and \cite[Section~6]{LiuYang2017} that the canonically induced metric $\omega$ on the diagonal Hopf surface $M=\S^1\times \S^3$ is Gauduchon and  $\mathrm{\Theta^{(2)}}(\omega)=\frac{1}{4}\omega>0$. However,   the top intersection number $c_1^2(M)=0$ and $H^2_{\mathrm{dR}}(M,\R)=H^{1,1}_{\bp}(M,\C)=0$. \\

	\noindent  On the other hand, Theorem \ref{secondChernRicci} is related to a question raised by Angella, Calamai and Spotti  in \cite[Remark~13]{Angella2020}. Actually, all weakly Hermitian-Einstein  surfaces are conformal K\"ahler except Hopf surfaces (\cite[Theorem~2]{Gauduchon1997}).  One can see that Hopf surfaces can not support weakly Hermitian-Einstein mertic with  Hermitian-Einstein function $u\leq 0$ (Theorem \ref{generalHE}).  More generally,

	\bcorollary\label{Hopf} Let $M$ be a Hopf surface. Then it can not support a smooth Hermitian metric $\omega$ with non-positive second Chern-Ricci curvature  $\mathrm{\Theta^{(2)}}(\omega)$.
	\ecorollary

	\noindent For more discussions on Chern scalar curvature and curvature relations,  we refer to \cite{Angella2017},  \cite{Chen2021}, 
	\cite{Chen2022},	\cite{Lejmi2018} and the references therein.
	Moreover,  by using the techniques in   \cite{Apostolov1999}, \cite{Chen2023}, \cite{Fu2022},  \cite{Li2010} and \cite{Tosatti2007}, one can  extend these results to almost Hermitian manifolds. In a subsequent paper \cite{Yang2024}, we shall discuss the complex geometry of various Riemannian curvatur tensors in details. \\

	

	\vskip 2\baselineskip

		\section{Background materials}
		In this section, we fix the conventions briefly for readers' convenience and also refer to \cite{LiuYang2017} and \cite{Yang2020MZ} for more details.
	Let $(M, g, \nabla)$ be a $2n$-dimensional Riemannian manifold with
	 Levi-Civita connection $ \nabla$. The tangent bundle of $M$ is denoted by $T_\R M$. The Riemannian curvature tensor of
	$(M,g, \nabla)$ is defined as $$
	R(X,Y,Z,W)=g\left( \nabla_X \nabla_YZ- \nabla_Y \nabla_XZ- \nabla_{[X,Y]}Z,W\right)$$
	for tangent vectors $X,Y,Z,W\in \Gamma(M,T_\R M)$. Let $T_\C M=T_\R M\ts \C$. One can extend the metric $g$ and the
	Levi-Civita connection $ \nabla$ to $T_{\C}M$ in the $\C$-linear way.
	Let $(M,g,J)$ be a Hermitian manifold and $z^i=x^i+\sq x^I$ be  local holomorphic coordinates on $M$.
	\noindent There is a Hermitian form $h:T_\C M\times T_\C M\>\C$ given
	by $ h(X,Y):= g(X,Y)$. The
	fundamental $2$-form associated to the $J$-invariant metric $g$ is  $\omega=\sq h_{i\bar j} dz^i\wedge d\bar z^j$. The complexified Christoffel symbols are given by \beq
	\Gamma_{AB}^C=\sum_{E}\frac{1}{2}g^{CE}\big(\frac{\p g_{AE}}{\p
		z^B}+\frac{\p g_{BE}}{\p z^A}-\frac{\p g_{AB}}{\p
		z^E}\big)=\sum_{E}\frac{1}{2}h^{CE}\big(\frac{\p h_{AE}}{\p
		z^B}+\frac{\p h_{BE}}{\p z^A}-\frac{\p h_{AB}}{\p z^E}\big)\label{realcomplexchristoffelsymbol}
	\eeq
	where $A,B,C,E\in \{1,\cdots,n,\bar{1},\cdots,\bar{n}\}$ and
	$z^{A}=z^{i}$ if $A=i$, $z^{A}=\bar{z}^{i}$ if $A=\bar{i}$. For
	example \beq \Gamma_{ij}^k=\frac{1}{2}h^{k\bar \ell}\left(\frac{\p
		h_{j\bar \ell}}{\p z^i}+\frac{\p h_{i\bar \ell}}{\p z^j}\right),\
	\Gamma_{\bar ij}^k=\frac{1}{2}h^{k\bar \ell}\left(\frac{\p h_{j\bar
			\ell}}{\p \bar z^i}-\frac{\p h_{j\bar i}}{\p \bar
		z^\ell}\right).\label{christoffelsymbols} \eeq Since $\Gamma_{AB}^C=\Gamma_{BA}^C$, we have  $\Gamma_{\bar i j}^k=\Gamma_{j\bar i}^k$. We also have $\Gamma_{\bar i\bar
		j}^k=\Gamma_{ij}^{\bar k}=0$ since
	$h_{pq}=h_{\bar i\bar j}=0$.
	

	\vskip 1\baselineskip

	\subsection{The induced Levi-Civita connection $\nabla^{\mathrm{LC}}$ on $(T^{1,0}M,h)$}  
	There is an induced connection
	$\nabla^{\mathrm{LC}}$ on  $T^{1,0}M$
	given by \beq \nabla^{\text{LC}}=\pi\circ\nabla:
	\Gamma(M,T^{1,0}M)\stackrel{\nabla}{\rightarrow}\Gamma(M, T_{\C}M\ts
	T_{\C}M)\stackrel{\pi}{\rightarrow}\Gamma(M,T_{\C}M\ts T^{1,0}M).
	\eeq  It is easy to see that
	 $\nabla^{\text{LC}}$ is
	determined by the following relations \beq
	\nabla^{\text{LC}}_{\frac{\p}{\p z^i}}\frac{\p}{\p
		z^k}:=\Gamma_{ik}^p\frac{\p}{\p z^p} \qtq{and}
	\nabla^{\mathrm{LC}}_{\frac{\p}{\p \bar z^j}}\frac{\p}{\p
		z^k}:=\Gamma_{\bar jk}^p\frac{\p}{\p z^p} \eeq where the Christoffel symbols are defined in (\ref{christoffelsymbols}). As usual, the curvature  tensor $\mathfrak{R}\in
	\Gamma(M,\Lambda^2 T_{\C}M\ts T^{*1,0}M\ts T^{1,0}M)$ of
	$\nabla^{\mathrm{LC}}$ is defined as \beq  \mathfrak{R}(X,Y)s
	=\nabla^{\mathrm{LC}}_{X}\nabla^{\mathrm{LC}}_Ys-\nabla^{\mathrm{LC}}_Y\nabla^{\mathrm{LC}}_Xs-\nabla^{\mathrm{LC}}_{[X,Y]}s \label{curvaturedefinition} \eeq
	for any $X,Y\in \Gamma(M,T_{\C}M)$ and $s\in \Gamma(M, T^{1,0}M)$. A straightforward
	computation shows that the curvature tensor $\mathfrak{R}$ has
	$(1,1)$ components \beq \mathfrak{R}_{i\bar
		jk}^{\ell}=-\left(\frac{\p \Gamma^{\ell}_{ik}}{\p \bar z^j}-\frac{\p
		\Gamma^{\ell}_{\bar jk}}{\p z^i}+\Gamma_{ ik}^{s}\Gamma^{\ell}_{\bar
		js}-\Gamma_{ \bar jk}^{s}\Gamma^{\ell}_{is }\right).\label{levicivitacurvatureformula}\eeq
	
	\vskip 1\baselineskip 	
	
	\subsection{The Chern connection $\nabla^{\mathrm{ch}}$ on $(T^{1,0}M,h)$} 
	\noindent 	The Chern
	connection $\nabla^{\mathrm{ch}}$ on $(T^{1,0}M,h)$ is  determined by the rules
	\beq
	\nabla^{\mathrm{ch}}_{\frac{\p}{\p z^i}}\frac{\p}{\p
		z^k}:={^c\Gamma}_{ik}^p\frac{\p}{\p z^p} \qtq{and}
	\nabla^{\mathrm{ch}}_{\frac{\p}{\p \bar z^j}}\frac{\p}{\p
		z^k}:=0 \eeq
	where ${^c\Gamma}_{ik}^p=h^{p\bar\ell}\frac{\p h_{k\bar\ell}}{\p z^i}$. 
	The curvature  components of $\nabla^{\mathrm{ch}}$ are \beq \Theta_{i\bar j k\bar
		\ell}=-\frac{\p^2 h_{k\bar \ell}}{\p z^i\p \bar z^j}+h^{p\bar
		q}\frac{\p h_{p\bar \ell}}{\p \bar z^j}\frac{\p h_{k\bar q}}{\p
		z^i}.\label{cherncurvatureformula} \eeq
	It is
	well-known that the \emph{first Chern-Ricci curvature} is \beq
	\Theta^{(1)}:= \sq\Theta^{(1)}_{i\bar j} dz^i\wedge d\bar z^j,\ \ \
	\Theta^{(1)}_{i\bar j}= h^{k\bar \ell}\Theta_{i\bar j k\bar \ell}
	=-\frac{\p^2 \log \det(h_{k\bar \ell})}{\p z^i\p\bar
		z^j}.\label{firstchernriccicurvatureformula}\eeq The \emph{second Chern-Ricci curvature}
	$\Theta^{(2)}=\sq \Theta^{(2)}_{i\bar j}dz^i\wedge d\bar z^j$ has
	components $$ \Theta^{(2)}_{i\bar j}=h^{k\bar \ell}\Theta_{k\bar
		\ell i\bar j}.$$

	\vskip 2\baselineskip

	\section{Weakly Hermitian-Einstein  surfaces}
	
	In this section, we refine the classical result of Gauduchon-Ivanov \cite{Gauduchon1997}.  On a compact Hermitian surface $\left(M,\omega=\sq
	h_{i\bar j}dz^i\wedge d\bar z^j\right)$, the torsion tensor $T$ of  $(M,h)$ has components
	\beq T_{ij}^k=h^{k\bar \ell}\left( \frac{\p h_{j\bar\ell}}{\p z^i}-\frac{\p h_{i\bar\ell}}{\p z^j}\right). \eeq
	We aslo write 
	\beq T_i:=\sum\limits_k T_{ik}^k,\ \ \ T_{\bar i}:=\bar {T_i}. \eeq 
	\noindent	
	It is easy to see that 
	\beq  \Gamma_{\bar i j}^k=\frac{1}{2}\bar{T_{ip}^q}h_{q\bar k}h^{j\bar p},\ \ \  2\sum_{k}\Gamma_{\bar i k}^k=T_{\bar i}.\eeq 	By the
well-known Bochner formula (e.g. \cite[Lemma~4.3]{LiuYang2012}), \beq [\bp^*,L]=\sq
\left(\p+\tau\right)\label{key0}\eeq  where $L \phi=\omega\wedge \phi$ and  $ \tau=[\Lambda,\p\omega]$,  one deduces that  \beq  \bp^*\omega=\sq T_{i}dz^i,\ \ \ \p^*\omega=-\sq T_{\bar j}d\bar z^j.\label{key1}
\eeq 
	Let $\theta:=T_idz^i+T_{\bar j}d\bar z^j$ be  the Lee form. By using \cite[Proposition ~2.2]{Tod1992},  Gauduchon and Ivanov obtained the following result in \cite[Theorem~1]{Gauduchon1997}.
\blemma  Let $(M,\omega)$ be a compact Hermitian surface. Suppose $\omega$ is Gauduchon and weakly Hermitian-Einstein, i.e.  $\p\bp\omega=0$ and
\beq \Theta^{(2)}(\omega)=u \omega \eeq 
for some $u\in C^\infty(M,\R)$. Then $\nabla\theta$ is skew-symmetric.
\elemma 

\noindent Actually,  it is proved in \cite[Theorem~1]{Gauduchon1997} that  when $\omega$ is Gauduchon and weakly Hermitian-Einstein, then the background Riemannian metric is Weyl-Einstein. By  \cite[Proposition ~2.2]{Tod1992}, they deduced that $\nabla\theta$ is skew-symmetric.  

\bcorollary Let $(M,\omega)$ be a Hermitian surface. If $\nabla \theta$ is skew-symmetric, then 
\beq \p\p^*\omega+\bp\bp^*\omega=-2\sq \bp^*\omega\wedge \p^*\omega+2|\p^*\omega|^2\omega. \label{S}\eeq 
In particular, $\omega$ is Gauduchon.
\ecorollary
\bproof  A straightforward computation shows that  
\be\left(\nabla \theta\right)\left(\frac{\p }{\p z^i},\frac{\p}{\p\bar z^j}\right) +\left(\nabla \theta\right)\left(\frac{\p}{\p\bar z^j},\frac{\p }{\p z^i}\right)&=&\nabla_i T_{\bar j}+\nabla_{\bar j} T_i-\Gamma_{i\bar j}^s T_s-\Gamma_{\bar j i}^{\bar s} T_{\bar s}\\&=&\frac{\p T_{\bar j}}{\p z^i}+\frac{\p T_{ i}}{\p\bar z^j}-2\Gamma_{i\bar j}^s T_s-2\Gamma_{\bar j i}^{\bar s} T_{\bar s}. \ee
Moreover, on a Hermitian surface, we have 
\beq 2\Gamma_{i\bar j}^s T_s=2\Gamma_{\bar j i}^{\bar s} T_{\bar s}=T_{i}T_{\bar j}-|\p^*\omega|^2 h_{i\bar j}. \eeq
Hence, we obtain  (\ref{S}).  By taking trace, we deduce that $$\Lambda(\p\p^*\omega+\bp\bp^*\omega)=2|\p^*\omega|^2.$$ Note also that $$\Lambda\p\p^*\omega-|\p^*\omega|^2=\sq \bp^*\p^*\omega.$$
Hence, $\omega$ is Gauduchon.
\eproof

\blemma Let $(M,\omega)$ be a compact Hermitian surface. Suppose $\omega$ is  Gauduchon and weakly Hermitian-Einstein. Then 
\beq \|\bp\bp^*\omega\|^2=\|\Lambda\bp\bp^*\omega\|^2=(|\p^*\omega|^4,1).\eeq 
\elemma 

\bproof  A straightforward calculation shows that 
\beq \left(\p\p^*\omega, \bp\bp^*\omega \right)=\|\Lambda\p\p^*\omega\|^2.\eeq 
Indeed, for any $(1,1)$ form $\eta$ on a Hermitian surface,  we have 
\beq *\eta=-\eta+(\Lambda \eta)\omega \eeq
and  so
\be 
\left(\p\p^*\omega, \bp\bp^*\omega \right)&=&\left(\partial^{*} \omega, \partial^{*} \overline{\partial \partial}^{*} \omega\right)=\left(\partial^{*} \omega,-* \bar{\partial} *\left(\overline{\partial \partial}^{*} \omega\right)\right) \\
&=&\left(\partial^{*} \omega, *\bp \bp \bp^{*} \omega\right)+\left(\partial^{*} \omega,-* \bar{\partial}\left(\left(\Lambda \bp \bp^{*} \omega\right) \omega\right)\right) \\
&=&\left(\bar{\partial}^{*} *\partial^{*} \omega,\left(\Lambda \bp \bp^{*} \omega\right) \omega\right)=\left(* \partial \partial^{*} \omega,\left(\Lambda \bp \bp^{*}\omega\right) \omega\right) \\
&=&\left(\Lambda \p \p^{*}\omega, \Lambda  \bp \bp^{*}\omega\right) .
\ee 
Since $\omega$ is Gauduchon, 
\beq \Lambda \p \p^{*}\omega=\Lambda  \bp \bp^{*}\omega=|\bp^*\omega|^2. \eeq 
Therefore, we obtain
\beq \left(\p\p^*\omega+\bp\bp^*\omega,\p\p^*\omega+\bp\bp^*\omega \right)=2\|\p\p^*\omega\|^2+2\|\Lambda\p\p^*\omega\|^2. \eeq 
On the other hand, by (\ref{S}), $\p\p^*\omega+\bp\bp^*\omega=2|\p^*\omega|^2\omega-2\sq \bp^*\omega\wedge \p^*\omega$, one has
$$\|\p\p^*\omega+\bp\bp^*\omega\|^2=\|2|\p^*\omega|^2\omega-2\sq \bp^*\omega\wedge \p^*\omega\|^2=4(|\p^*\omega|^4,1)$$
and so $\|\bp\bp^*\omega\|^2=\|\Lambda\bp\bp^*\omega\|^2=(|\p^*\omega|^4,1)$.
\eproof 

\noindent The following result is proved in \cite[Theorem~2]{Gauduchon1997}.  We give a simplified proof here.
	
	\btheorem\label{GI} Let $(M,\omega)$ be a compact Hermitian surface. Suppose $\omega$ is  Gauduchon and weakly Hermitian-Einstein, i.e. 
	\beq \Theta^{(2)}(\omega)=u \omega \eeq 
	for some $u\in C^\infty(M,\R)$. Then $u-|\bp^*\omega|^2$ is a constant. Moreover, either  $\omega$ is a K\"ahler metric, or  $u=|\bp^*\omega|^2$.
	\etheorem
	
	\bproof In this case, we have 
	\beq \Theta^{(1)}-\left(\partial \partial^{*} \omega+\overline{\partial \partial}^{*} \omega\right)=\left( u-|\bp^*\omega|^2\right)\omega. \eeq 
	Let $f=u-|\bp^*\omega|^2$. By applying $\p\bp$ to both sides, we have 
	\beq 0=\p\bp(f\omega)=\p\bp f\wedge \omega+\p f\wedge \bp\omega-\bp f\wedge \p\omega. \eeq 
	By taking trace, we get an elliptic equation. The maximum principle tells that $f$ is a constant. ( One can also deduce that by multiplying $f$ and doing integration by parts.) We denote this constant by $\lambda$ and so 
	$$\Theta^{(1)}-\left(\partial \partial^{*} \omega+\overline{\partial \partial}^{*} \omega\right)=\lambda \omega.$$
By taking $\p$ on both sides, we deduce that $-\p\bp\bp^*\omega=\lambda\p\omega$.
Therefore $$-\left(\partial \overline{\partial \partial}^{*} \omega, \partial \omega\right)=\lambda(\partial \omega, \partial \omega)=\lambda\left\|\partial^{*} \omega\right\|^2.$$ On the other hand,
\be
\left(\partial \overline{\partial \partial}^{*} \omega, \partial \omega\right)&=&\left(\overline{\partial \partial}^{*} \omega, \partial^{*} \partial \omega\right)=\left(\overline{\partial \partial}^{*} \omega,-* \bar{\partial} * \partial * \omega\right) \\
&=&\left(\overline{\partial \partial}^{*} \omega, * \overline{\partial \partial^{*}} \omega\right)=-\|\overline{\partial \partial^{*}} \omega\|^2+\|\Lambda\bp\bp^*\omega\|^2 \\
&=&-\|\overline{\partial \partial^{*}} \omega\|^2+\left(|\overline{\partial}^{*} \omega|^4, 1\right)=0
\ee 
since  $\|\bp\bp^*\omega\|^2=\|\Lambda\bp\bp^*\omega\|^2=\left(|\overline{\partial}^{*} \omega|^4, 1\right)$. Therefore, we obtain
	\beq\lambda \|\bp^*\omega\|^2=0.\eeq 
If $\lambda\neq 0$, then $\bp^*\omega=0$ and so $\omega$ is K\"ahler. If $\lambda=0$, then 	$u-|\bp^*\omega|^2=0$.
	\eproof

\noindent 	Actually, it is proved in \cite{Gauduchon1995} that if $u-|\bp^*\omega|^2=0$, then $M$ is a Hopf surface which is non-K\"ahler. In this case 
\beq  \Theta^{(1)}-\left(\partial \partial^{*} \omega+\overline{\partial \partial}^{*} \omega\right)=0.\eeq 
We refine Theorem \ref{GI} and obtain a generalization of Theorem \ref{main0}:

\bproposition\label{generalGI}   Let $(M,\omega)$ be a compact Hermitian surface. Suppose $\omega$ is  Gauduchon and weakly Hermitian-Einstein,
\beq \Theta^{(2)}(\omega)=u \omega \eeq 
for some $u\in C^\infty(M,\R)$.  If there exists some Hermitian metric $\omega_0$ such that $\int_M u\omega_0^2\leq 0$, then $\omega$ is K\"ahler.
\eproposition
	
\bproof  As in the proof of Theorem \ref{GI}, $$u=\lambda+|\p^*\omega|^2$$ for some constant $\lambda$.  If there exists some Hermitian metric $\omega_0$ such that $\int_M u\omega_0^2\leq 0$, we deduce that either $\lambda=0$ or $\lambda<0$. If $\lambda=0$, we have $|\p^*\omega|^2=0$ and so $\omega$ is K\"ahler. If $\lambda<0$, by Theorem \ref{GI}, $\omega$ is K\"ahler.
\eproof 
	
	\noindent 
	The following results are esstential known to experts along this line.
	
	\btheorem\label{generalHE} Let $(M,\omega)$ be a compact Hermitian surface. Suppose $\omega$ is weakly Hermitian-Einstein, i.e. 
	\beq \Theta^{(2)}(\omega)=u \omega \eeq 
	for some $u\in C^\infty(M,\R)$. If $u\leq 0$, then $M$ is a K\"ahler surface and $\omega$ is  conformal to a K\"ahler metric.
	\etheorem
	
	\bproof  Let $\omega_f=e^f\omega$ be the Gauduchon metric in the conformal class of $\omega$.  By using the curvature formula of the Chern connection, one has  
	\be \Theta_{i\bar j k\bar
		\ell}(\omega_f)&=&-\frac{\p^2 (e^f h_{k\bar \ell})}{\p z^i\p \bar z^j}+e^{-f}h^{p\bar
		q}\frac{\p(e^f h_{p\bar \ell})}{\p \bar z^j}\frac{\p (e^fh_{k\bar q})}{\p
		z^i}=e^f\Theta_{i\bar j k\bar\ell }-e^f h_{k\bar \ell } \frac{\p^2 f}{\p z^i\p \bar z^j}.  \ee  In particular,
	\beq \Theta^{(2)}(\omega_f)=\Theta^{(2)}(\omega)-\left(\sqrt{-1} \Lambda_\omega \partial \bar{\partial} f\right) \omega=e^{-f}\left(u-\left(\sqrt{-1} \Lambda_\omega \partial \bar{\partial} f\right)\right) \omega_f. \eeq 
For simplicity, we set \beq u_f:=e^{-f}u- \mathrm{tr}_{\omega_f} (\sq \p\bp f). \eeq 	
	It is obvious that $$\int_M u_f \omega_f^2=\int_M e^{-f}u \omega_f^2-\int_M \mathrm{tr}_{\omega_f} (\sq \p\bp f) \omega_f^2=\int_M e^{-f}u \omega_f^2\leq 0.$$
	By Proposition \ref{generalGI}, $\omega_f=e^f\omega$ is K\"ahler.
	\eproof

		\bcorollary  Let $(M,\omega)$ be a compact Hermitian surface with $\Theta^{(2)}(\omega)=-c^2\omega$ for some $c\in \R$. 
	Then $(M,\omega)$ is  K\"ahler-Einstein.
	\ecorollary 
	
	\bproof  Suppose  that $\omega_f=e^f\omega$ is a Gauduchon metric in the conformal class of $\omega$ with $\inf_M f=0$. By Theorem \ref{generalHE}, we deduce that $$\Theta^{(2)}(\omega_f)=u_f\omega_f$$ and $\omega_f$ is K\"ahler.  Since $\Theta^{(1)}(\omega_f)=\Theta^{(2)}(\omega_f)$ is $\p\bp$-closed, we have   $$\mathrm{tr}_{\omega_f} \sq \p\bp u_f=0.$$ The strong maximum principle asserts that $u_f$ is a constant. Since $f\geq 0$ and $$u_f=-e^{-f}c^2- \mathrm{tr}_{\omega_f} (\sq \p\bp f),\ \  \int_M u_f \omega_f^2=-\int_M e^{-f}c^2 \omega_f^2\leq 0,$$ we deduce that $$-c^2\leq u_f\leq 0.$$ On the other hand, at a minimum point $p$ of $f$, $f(p)=0$, and we have 
	$$u_f=-c^2- \mathrm{tr}_{\omega_f} (\sq \p\bp f)(p)\leq -c^2.$$
	Hence, $u_f=-c^2$. Therefore,  $f=0$ and $\omega_f=\omega$ is K\"ahler. 
	\eproof

		\vskip 2\baselineskip

			\section{Manifolds with non-positive second Chern-Ricci curvature}
			
			In this section, we prove Theorem \ref{secondChernRicci}, Corollary \ref{surface} and Corollary \ref{Hopf}. The key ingredient is the following theorem, which is obtained by
		 using the Chern-Weil theory for a combination of metric connections on  Hermitian line bundles.
			
		\btheorem\label{thm:ChernWeil} Let $M$ be a compact complex manifold with $\dim_\C M=n$ and $L$ be a holomorphic line bundle over $M$. For $i=1,\cdots, m$,  let $h_i$ be a smooth Hermitian metric on $L$ and $\nabla_i$ be a metric connection on $(L,h_i)$.  Let $R_i$ be  the $(1,1)$-component  of the curvature tensor of $(L,h_i,\nabla_i)$. Suppose that there exist real numbers $a_1,\cdots, a_m$ such that 
		$$\sum_{i=1}^m a_i R_{i}$$
		is a quasi-positive $(1,1)$-form over $M$.
		\bd \item If $\sum_{i}a_i>0$, then the top interesection number $c^n_1(L)>0$;
		
		\item If $\sum_{i} a_i<0$, then the signed  top interesection number $(-1)^nc^n_1(L)>0$.
		\ed
		\etheorem

		\bproof Let $\nabla$ be an arbitrary metric compatible connection on the Hermitian holomorphic line bundle $(L,h)$ and $R\in \Gamma(M,\Lambda^2T^*M)$ be its curvature tensor. There is a natural decomposition according to the degrees
		\beq  R=\eta^{2,0}+\eta^{0,2}+\eta^{1,1}.\eeq
		Since $\nabla$ is metric compatible, we know $$\eta^{0,2}=\bar{\eta^{2,0}}.$$
		Let $\nabla^{\mathrm{ch}}$ be the Chern connection on $(L,h)$ and $\Theta\in\Gamma(M,\Lambda^{1,1}T^*M)$ be its Chern curvature. By the Chern-Weil theory (e.g. \cite{Zhang2001}), one has 
		\beq R-\Theta=d\beta \eeq
		for some $1$-form $\beta$ on $M$.  We apply this to the family $\{(L,h_i,\nabla_i)\}_{i=1}^m$ and deduce that 
		$$\sum a_i R_i- \sum a_i \Theta=d\beta_0+\eta_0^{2,0}+\bar{\eta_0^{2,0}}$$
		for some $1$-form $\beta_0$ and $(2,0)$-form $\eta^{2,0}_0$. Note that both of them depend on $a_i$.\\
		
		Suppose $a=\sum_i a_i\neq 0$. We consider a $2$-form
		\beq W:=\frac{1}{a}\left(\sum_i a_i R_i-\eta_0^{2,0}-\bar{\eta_0^{2,0}}\right).\eeq 
		It is clear that \beq  W-\Theta=a^{-1}d\beta_0, \eeq and 
		\be c_1^n(L)=\int_M \Theta^n=\int_M W^n
		=
		a^{-n}\sum_{\ell=0}^{\left[\frac{n}{2}\right]}{n\choose
			{2\ell}}{{2\ell}\choose \ell} \int_M(\eta_0^{2,0}\wedge \bar
		{\eta_0^{2,0}})^{\ell}\wedge \left(\sum_i a_i R_i\right)^{n-2\ell}. \ee
		If  $\suml_i a_i R_i$ is quasi-positive, then
		$$\sum_{\ell=0}^{\left[\frac{n}{2}\right]}{n\choose
			{2\ell}}{{2\ell}\choose \ell} \int_M(\eta_0^{2,0}\wedge \bar
		{\eta_0^{2,0}})^{\ell}\wedge \left(\sum_i a_i R_i\right)^{n-2\ell}\geq 0$$
		for $1\leq \ell\leq
		\left[\frac{n}{2}\right]$. Moreover, one has  $$\int_M  \left(\sum_i a_i R_i\right)^{n}>0 .$$
		Hence, we obtain $c_1^n(L)>0$ if $a>0$ and 	$(-1)^nc_1^n(L)>0$ if $a<0$.\eproof

		\noindent We need the following curvature estimate in the proof of Theorem \ref{secondChernRicci}.
		\blemma \label{pluriclosedcurvaturerelation} Let $(M,\omega)$ be a Hermitian manifold with $\Lambda_\omega \p\bp \omega=0$. Then
		\beq 2\ssR^{(1)}-\Theta^{(1)}\leq  \Theta^{(2)} \eeq
		as Hermitian $(1,1)$-forms on $M$, where 
		$$ \mathfrak{R}^{(1)}=\sq \mathfrak R_{i\bar j k}^kdz^i\wedge d\bar
		z^j.$$
		 Moreover, the identity holds if and only if $\omega$ is K\"ahler.
		\elemma
		
		\bproof By using the curvature formula, we have 
		$$\Theta^{(1)}_{i\bar j}=-\frac{\p {^c\Gamma_{ik}^k}}{\p \bar z^j},\ \ \Theta^{(2)}_{i\bar j}=-h_{p\bar j}h^{k\bar \ell}\frac{\p{^c\Gamma_{ki}^p}}{\p\bar z^\ell} \qtq{and} \ssR^{(1)}_{i\bar j}\stackrel{(\ref{levicivitacurvatureformula})}{=}-\left(\frac{\p\Gamma_{ik}^k}{\p\bar z^j}-\frac{\p\Gamma_{\bar j k}^k	}{\p z^i}\right).$$
		Hence, 
		\be &&\Theta^{(1)}+\Theta^{(2)}-2\ssR^{(1)}=\frac{\p {^c\Gamma}_{ki}^k}{\p \bar z^j}-2\frac{\p\Gamma_{\bar j k}^k	}{\p z^i} -h_{p\bar j}h^{k\bar \ell}\frac{\p{^c\Gamma_{ki}^p}}{\p\bar z^\ell}\\
		&=&\frac{\p }{\p \bar z^j}\left(h^{k\bar m}\frac{\p h_{i\bar m}}{\p z^k}\right)-\frac{\p}{\p  z^i}\left(h^{k\bar m}\left(\frac{\p h_{k\bar m}}{\p \bar z^j}-\frac{\p h_{k\bar j}}{\p \bar z^m}\right) \right) -h_{p\bar j}h^{k\bar \ell}\frac{\p}{\p\bar z^\ell}\left(h^{p\bar m}\frac{\p h_{i\bar m}}{\p z^k}\right)\\&\stackrel{(*)}{=}&\frac{\p h^{k\bar m} }{\p \bar z^j}\cdot \frac{\p h_{i\bar m}}{\p z^k}-\frac{\p h^{k\bar m}}{\p  z^i}\left(\frac{\p h_{k\bar m}}{\p \bar z^j}-\frac{\p h_{k\bar j}}{\p \bar z^m}\right) -h_{p\bar j}h^{k\bar \ell}\frac{\p h^{p\bar m}}{\p\bar z^\ell}\cdot \frac{\p h_{i\bar m}}{\p z^k}\\
		&\stackrel{(**)}{=}& h^{p\bar q}
		h_{k \bar \ell}\cdot T_{ip}^k\cdot \bar {T_{jq}^\ell},\ee	
		where we use the condition $\Lambda_\omega \p\bp \omega=0$ in the identity $(*)$. Indeed,
		\be &&h^{k\bar m}\frac{\p^2 h_{i\bar m}}{\p\bar z^j\p z^k}- h^{k\bar m}\frac{\p}{\p  z^i}\left(\frac{\p h_{k\bar m}}{\p \bar z^j}-\frac{\p h_{k\bar j}}{\p \bar z^m}\right) -h_{p\bar j}h^{k\bar \ell} h^{p\bar m} \frac{\p^2 h_{i\bar m}}{\p\bar z^\ell\p z^k}\\
		&=&h^{k\bar m}\left(\frac{\p^2 h_{i\bar m}}{\p\bar z^j\p z^k}+\frac{\p^2 h_{k\bar j}}{\p\bar z^m\p z^i}-\frac{\p^2 h_{k\bar m}}{\p\bar z^j\p z^i}-\frac{\p^2 h_{i\bar j}}{\p\bar z^m\p z^k}\right)=(\Lambda_\omega\p\bp\omega)_{i\bar j}.\ee
		The verification of $(**)$ is straightforward.  Hence, we have 
		\beq \Theta^{(1)}+\Theta^{(2)}=2\ssR^{(1)}+Q\eeq
		where $Q:= \sq Q_{i\bar j}dz^i\wedge d\bar z^j=\sq h^{p\bar q}
		h_{k \bar \ell}\cdot T_{ip}^k\cdot \bar {T_{jq}^\ell} dz^i\wedge d\bar
		z^j.$  It is obvious that $Q$ is a Hermitian semi-positive $(1,1)$-form on $M$. Moreover, if  $ \Theta^{(1)}+\Theta^{(2)}=2\ssR^{(1)}$, then $Q\equiv 0$ and $T\equiv 0$. Hence, $\omega$ is K\"ahler.
		\eproof

		\noindent \emph{Proof of Theorem \ref{secondChernRicci}.}	We apply the Chern-Weil theory to the Hermitian manifold
		$(M,\omega)$. Let $E=T^{1,0}M$ and $h$ be the Hermitian metric  induced by	$\omega$. With respect to the Levi-Civita connection $\nabla^{\mathrm{LC}}$
		on $(E,h)$, we have a decomposition $$ R^{\det
			E}=\eta^{2,0}+\eta^{0,2}+\eta^{1,1}.$$  It is obvious that
		\beq \eta^{1,1}= \sq \mathfrak R_{i\bar j k}^kdz^i\wedge d\bar
		z^j=\mathfrak{R}^{(1)},\ \ \
		\qtq{and}
		\eta^{0,2}=\bar{\eta^{2,0}}.\eeq
		If $\Lambda\p\bp\omega=0$ and $\Theta^{(2)}$ is quasi-positive, by Lemma \ref{pluriclosedcurvaturerelation},
		$$-2\ssR^{(1)}+\Theta^{(1)}\geq -\Theta^{(2)}$$
		is quasi-positive. By Theorem \ref{thm:ChernWeil}, we obtain the conclusion that $(-1)^n c_1^n(M)>0$. This implies the non-vanishing of the cohomology groups $H^2_{\mathrm{dR}}(M,\R)$, $H^{1,1}_{\bp}(M,\C)$ and $H^{1,1}_{\mathrm{BC}}(M,\C)$. \qed

		\vskip 1\baselineskip

		\noindent \emph{Proof of Corollary \ref{surface}.} By Lemma \ref{pluriclosedcurvaturerelation},
		$$-2\ssR^{(1)}+\Theta^{(1)}= -\Theta^{(2)}+Q.$$
		On a Hermitian surface, we have 
		\beq Q=\sq h^{p\bar q}
		h_{k \bar \ell}\cdot T_{ip}^k\cdot \bar {T_{jq}^\ell} dz^i\wedge d\bar
		z^j=\left|\partial^{*} \omega\right|^2 \omega\eeq 
		 Indeed, we can verify it by using normal coordinates (e.g. \cite[Lemma~3.4]{LiuYang2017}) at a point $p\in M$. That is,  there exist local holomorphic ``normal
		coordinates" $\{z^i\}$ centered at $p$ such that \beq h_{i\bar
			j}(p)=\delta_{ij},\ \ \ \Gamma_{ij}^k(p)=0,\ \ \ \Gamma_{\bar j
			i}^k(p)=\frac{\p h_{i\bar k}}{\p \bar z^j}(p)=-\frac{\p h_{i\bar
				j}}{\p\bar z^k}(p).\eeq
		If we set 
		\beq a_{1\bar 1}=\frac{\p h_{2\bar 2}}{\p z^1}\frac{\p h_{2\bar 2}}{\p \bar z^1},\ \ \ \ \ a_{1\bar 2}=\frac{\p h_{2\bar 2}}{\p z^1}\frac{\p h_{1\bar 1}}{\p \bar z^2},\ \ \ \ a_{2\bar 1}=\frac{\p h_{1\bar 1}}{\p z^2}\frac{\p h_{2\bar 2}}{\p \bar z^1},\ \ \ \ a_{2\bar 2}=\frac{\p h_{1\bar 1}}{\p z^2}\frac{\p h_{1\bar 1}}{\p \bar z^2}, \eeq 	
		one can compute easily that 	
		\beq |\bp^*\omega|^2=4\left(\frac{\p h_{2\bar 2}}{\p z^1}\frac{\p h_{2\bar 2}}{\p \bar z^1}+\frac{\p h_{1\bar 1}}{\p z^2}\frac{\p h_{1\bar 1}}{\p \bar z^2}\right) =4(a_{1\bar 1}+a_{2\bar 2})\eeq 
		and 
		\beq Q=\sum_{i,j,p, q}\sqrt{-1}  T_{i p}^q \overline{T_{j p}^q} d z^i \wedge d \bar{z}^j=4\sq (a_{1\bar 1}+a_{2\bar 2})(dz^1\wedge d\bar z^1+dz^2\wedge d\bar z^2) \eeq
		and so $Q=\left|\partial^{*} \omega\right|^2 \omega$. \\
		
		 From the proof of Theorem \ref{secondChernRicci}, we know 
		 \beq c^2_1(M)\geq \int_M \left( -\Theta^{(2)}+Q\right)^2\geq \int_M Q^2= \int_M |\p^*\omega|^4\omega^2\geq 0. \eeq 
		Hence, if $c_1^2(M)=0$, then $\p^*\omega=0$ and $\omega$ is K\"ahler.  Otherwise, $c^2_1(M)>0$ and by \cite[Theorem~9]{Kodaira1964}, we know $M$ is projective. \qed

		\vskip 1\baselineskip

	\noindent \emph{Proof of Corollary \ref{Hopf}.} Suppose $M$ has a Hermitian metric with $\Theta^{(2)}(\omega)\leq 0$.  Then the Chern scalar curvature $s_{\mathrm{\C}}\leq 0$. If $s_{\mathrm{\C}}= 0$, then $\Theta^{(2)}(\omega)=0$. By Theorem \ref{generalHE}, $M$ is K\"ahler and it is impossible. If $s_{\mathrm{\C}}$ is quasi-negative, then by \cite[Theorem~1]{Yang2019TAMS}, $K_M^{-1}$ is not pseudo-effective. But it is well-known that $K_{M}^{-1}$ is semi-positive. This is a contradiction. \qed

	\def\cprime{$'$} 

\end{document}